\documentclass[12pt]{amsart}
\usepackage{amsmath, amssymb,graphics,multirow}
\usepackage{hyperref}
\hypersetup{%
 colorlinks=true,
 linkcolor=blue,
}

\newtheorem{theorem}{Theorem}[section]
\newtheorem{lemma}[theorem]{Lemma}

\theoremstyle{definition}
\newtheorem{definition}[theorem]{Definition}

\newtheorem{remark}[theorem]{Remark}

\numberwithin{equation}{section}
\numberwithin{table}{section}

\newcommand{\bbQ}{{\mathbb Q}}
\newcommand{\bbR}{{\mathbb R}}
\newcommand{\bbZ}{{\mathbb Z}}

\newcommand{\cM}{{\mathcal M}}

\newcommand{\cO}{{\mathcal O}}

\def\Char{\operatorname{char}}

\begin{document}

\baselineskip=17pt

\title[Greatest common valuation of $\phi_{n}$ and $\psi_{n}^{2}$]
{The greatest common valuation of $\phi_{n}$ and $\psi_{n}^{2}$ at points on elliptic curves}

\author{Paul Voutier}
\address{}
\curraddr{London, UK}
\email{paul.voutier@gmail.com}

\author{Minoru Yabuta}
\address{}
\curraddr{Senri High School, 17-1, 2 chome, Takanodai, Suita, Osaka, 565-0861, Japan}
\email{yabutam@senri.osaka-c.ed.jp, rinri216@msf.biglobe.ne.jp}


\date{}

\dedicatory{}


\begin{abstract}
Given a minimal model of an elliptic curve, $E/K$, over a finite extension, $K$, of
$\bbQ_{p}$ for any rational prime, $p$, and any point $P \in E(K)$ of infinite
order, we determine precisely $\min \left( v \left( \phi_{n}(P) \right), v \left( \psi_{n}^{2}(P) \right) \right)$,
where $v$ is a normalised valuation on $K$ and $\phi_{n}(P)$ and $\psi_{n}(P)$
are polynomials arising from multiplication by $n$ on this model of the curve.
\end{abstract}

\maketitle

\section{Introduction}

Let $E/K$ be an elliptic curve given by a Weierstrass equation and let
$\left\{ \psi_{n}(X,Y) \right\}_{n \geq 1}$ be the sequence of associated
division polynomials. We can also define a related family of polynomials,
$\phi_{n}(X,Y)$, such that
\[
x([n]P) = \frac{\phi_{n}(x(P), y(P))}{\psi_{n}^{2}(x(P), y(P))},
\]
for any point, $P=(x(P),y(P)) \in E(K)$, satisfying the Weierstrass equation.

For many problems related to integral (and more generally, $S$-integral) points
on elliptic curves, it is important to know, or at least bound,
$\gcd \left( \phi_{n}(P), \psi_{n}^{2}(P) \right)$. For example, Ayad \cite{Ayad}
used such information to find all $S$-integral points on some elliptic curves of
rank $1$ (see his Sections~6--8 for three specific examples demonstrating his
technique). More recently, Stange \cite{Stange}, building on work of Ingram
\cite{Ingram2}, showed that if $E/\bbQ$ is an elliptic curve in minimal Weierstrass
form and $P \in E(\bbQ)$ is a non-torsion point, then there is at most one
value of $n$ larger than a bound that can be made effective such that $[n]P$
is integral. This leads to results (see Corollary~1.2 and the discussion that
follows in \cite{Stange}) related to the Hall-Lang Conjecture \cite{Lang}.

Good knowledge of $\gcd \left( \phi_{n}(P), \psi_{n}^{2}(P) \right)$ is also
required for problems involving elliptic divisibility sequences -- this is not
unrelated to the fact that the above two applications also involve rank
$1$ subgroups of $E(\bbQ)$. An example of this is Lemma~3.5 (and Remark~3.6) in
\cite{VY2}. A weaker version of that Lemma~3.5 would have led to a weaker version
of the main result there. In fact, trying to understand and generalise that Lemma~3.5
is how our interest in this subject arose.
In this paper, we determine this gcd precisely. This is a special case of the
more general result we prove here.

Our results are significantly better than previous results. Previous results
typically depend on the resultant of $\phi_{n}(P)$ and $\psi_{n}^{2}(P)$,
resulting in exponents that are $O \left( n^{4} \right)$ (e.g., Lemmas~4 and 5,
along with the claims used in their proofs, in \cite{Ingram1}), whereas here
our results are precise and we show that the actual growth of the exponents is
$O \left( n^{2} \right)$.

Throughout this paper, we let $p$ be a rational prime, $K$ a finite extension of $\bbQ_{p}$, $R$ the
ring of integers of $K$, with maximal ideal $\cM$, $\pi$ a uniformiser for
$R$ (i.e., $\cM = \pi R$), residue field $k=R/\cM$ and $v$ a valuation
for $K$ normalised so that $v(\pi) = 1$. 

We let $\lambda_{v}: E \left( K \right) \backslash \{ O \} \rightarrow \bbR$ be
the local height function for $E$ at $v$ as defined by Silverman in \cite{Silv4}.
We use this local height function to be consistent with the one we used in our
earlier papers. See \cite[Section~4]{CPS} for
helpful notes about the various normalisations of local height functions.

\begin{theorem}
\label{thm:main}
Assume we have a minimal Weierstrass model of an elliptic curve $E/K$ and that
$P \in E(K)$ is of infinite order. Let $n$ be a positive integer and put
\[
k_{v,n}(P) = \min \left( v\left( \phi_{n}(P) \right), v\left( \psi_{n}^{2}(P) \right) \right).
\]

If $P$ modulo $\pi$ is non-singular, then
$k_{v,n}(P)=\min \left( 0, n^{2}v(x(P)) \right)$.

If $P$ modulo $\pi$ is singular, then
\[
k_{v,n}(P) = -\frac{2[K:\bbQ_{p}]}{\log |k|}\lambda_{v}(P) n^{2}+\epsilon_{v,n}(P),
\]
where $\epsilon_{v,n}(P)$ is as in Table~$\ref{table:eps}$.

\begin{table}[h]
\centering
\scalebox{0.9}{%
\hspace{-10.0mm}\begin{tabular}{|c|c|c|cl|}
\hline
{\rm Kodaira symbol}          & $m_{P}$ & $-2[K:\bbQ_{p}]\lambda_{v}(P)/\log |k|$ & \multicolumn{2}{c|}{$\epsilon_{v,n}(P)$} \\ \hline
$III^{*}$                     &   $2$   &     $3/2$ & $0$        & \text{{\rm if} $n \equiv 0 \bmod{m_{P}}$}     \\
                              &         &           & $-3/2$     & \text{{\rm if} $n \not\equiv 0 \bmod{m_{P}}$} \\ \hline
$IV^{*}$                      &   $3$   &     $4/3$ & $0$        & \text{{\rm if} $n \equiv 0 \bmod{m_{P}}$}     \\
                              &         &           & $-4/3$     & \text{{\rm if} $n \not\equiv 0 \bmod{m_{P}}$} \\ \hline
$III$                         &   $2$   &     $1/2$ & $0$        & \text{{\rm if} $n \equiv 0 \bmod{m_{P}}$}     \\
                              &         &           & $-1/2$     & \text{{\rm if} $n \not\equiv 0 \bmod{m_{P}}$} \\ \hline
$IV$                          &   $3$   &     $2/3$ & $0$        & \text{{\rm if} $n \equiv 0 \bmod{m_{P}}$}     \\
                              &         &           & $-2/3$     & \text{{\rm if} $n \not\equiv 0 \bmod{m_{P}}$} \\ \hline
$I_{m}^{*}$, $c_{v}=2$        &   $2$   &     $1$   & $0$        & \text{{\rm if} $n \equiv 0 \bmod{m_{P}}$}     \\
                              &         &           & $-1$       & \text{{\rm if} $n \not\equiv 0 \bmod{m_{P}}$} \\ \hline
$I_{m}^{*}$, $m$ {\rm odd}, $c_{v}=4$
                              &   $2$   &     $1$   & $0$        & \text{{\rm if} $n \equiv 0 \bmod{m_{P}}$}     \\
$[2]P$ {\rm non-singular} $\bmod \, \pi$
                              &         &           & $-1$       & \text{{\rm if} $n \not\equiv 0 \bmod{m_{P}}$} \\ \hline
$I_{m}^{*}$, \,$m$ {\rm odd}, $c_{v}=4$,
                              &   $4$   & $(m+4)/4$ & $0$        & \text{{\rm if} $n \equiv 0 \bmod{m_{P}}$}   \\
$[2]P$ {\rm singular} $\bmod \, \pi$
                              &         &           & $-(m+4)/4$ & \text{{\rm if} $n \equiv 1,3 \bmod{m_{P}}$} \\
                              &         &           & $-1$       & \text{{\rm if} $n \equiv 2 \bmod{m_{P}}$}   \\ \hline
$I_{2m}^{*}$, $c_{v}=4$       &   $2$   & $v\left( \phi_{2}(P) \right)/4$ & $0$ & \text{{\rm if} $n \equiv 0 \bmod{m_{P}}$}   \\
                              &         &    & $-v\left( \phi_{2}(P) \right)/4$ & \text{{\rm if} $n \not\equiv 0 \bmod{m_{P}}$} \\ \hline
$I_{m}$                       & $\dfrac{m}{\gcd \left( a_{P,v}, m \right)}$
                              & $\dfrac{a_{P,v}\left( m-a_{P,v} \right)}{m}$
                              & $-\dfrac{n'\left( m-n' \right)}{m}$
                              & \text{{\rm if} $a_{P,v}n \equiv n' \bmod{m}$} \\ \hline
\end{tabular}%
}
\caption{$\epsilon_{v,n}(P)$ values}
\label{table:eps}
\end{table}
\end{theorem}

\begin{remark}
(i) $c_{v}$ is the size of the component group of $E$ at $v$, which we define
at the end of Subsection~\ref{subsect:notation}.

\vspace*{1.0mm}

(ii) $m_{P}$ is the smallest positive integer such that $\left[ m_{P} \right]P$
modulo $\pi$ is non-singular.

\vspace*{1.0mm}

(iii) If the Kodaira symbol is $I_{0}$, $I_{1}$, $II$ or $II^{*}$, then we cannot
have singular reduction. Hence these symbols do not appear in Table~\ref{table:eps}.

\vspace*{1.0mm}

(iv) In the entry for $I_{2m}^{*}$ with $c_{v}=4$ in Table~\ref{table:eps},
$v\left( \phi_{2}(P) \right)$ can take only two values, either $4$ or
$2m+4$ -- see Lemma~\ref{lem:ImStar-vals}.

\vspace*{1.0mm}

(v) $a_{P,v}$ in the entry for $I_{m}$ in Table~\ref{table:eps} is the
component of the N\'{e}ron model special fibre containing $P$. See Lemma~5.1 and
the surrounding text in \cite{Silv4} for more information.
We use $m$ instead of $m_{P}$ for the modulus here as it results in a
simpler expression for $\epsilon_{v,n}$.
We let $n'$ be the
smallest non-negative representative
of the congruence class $a_{P,v}n \bmod{m}$, so $0 \leq n'<m$.
\end{remark}

\begin{remark}
The referee has kindly shared with us other ways that the quantities in our
theorem arise. For example, when calculating height pairings for $P,Q \in E(K)$,
our $-\epsilon_{v,n}(P)$ is the correction term, $\mathrm{contr}_{v}(P)
=\mathrm{contr}_{v}(P,P)$, that arises in Section~11.8 of \cite{SS}. In the
notation of that section, we have $i=j$ here, so the values can be read from
the first row in Table~4 on the top of page~111 of \cite{SS}.
\end{remark}


\section{Preliminaries}


\subsection{Notation}
\label{subsect:notation}
Let $E/K$ be an elliptic curve given by the Weierstrass equation
\[
E/K: y^{2}+a_{1}x y+a_{3}y=x^{3}+a_{2}x^{2}+a_{4}x+a_{6},
\]
with $a_{1}, a_{2}, a_{3}, a_{4}, a_{6} \in R$.

We will also require the following quantities
\begin{align*}
b_{2}  &=  a_{1}^{2}+4 a_{2},   \\
b_{4}  &=  2 a_{4}+a_{1}a_{3},  \\
b_{6}  &=  a_{3}^{2}+4 a_{6},   \\
b_{8}  &=  a_{1}^{2}a_{6}+4 a_{2}a_{6}-a_{1}a_{3}a_{4}+a_{2}a_{3}^{2}-a_{4}^{2}, \\
c_{4}  &=  b_{2}^{2}-24b_{4},      \\
c_{6}  &= -b_{2}^{3}+36b_{2}b_{4}-216b_{6}, \\
\Delta &= -b_{2}^{2}b_{8}-8b_{4}^{3}-27b_{6}^{2}+9b_{2}b_{4}b_{6}, \\
j      &=  c_{4}^{3}/\Delta,
\end{align*}
where $\Delta$ is the {\it discriminant} of the Weierstrass equation. Note that
$4b_{8}=b_{2}b_{6}-b_{4}^{2}$ and $1728\Delta=c_{4}^{3}-c_{6}^{2}$.
If $\Char \left( \overline{K} \right) \neq 2$, then $E/K$ is also given by
$y^{2}=4x^{3}+b_{2}x^{2}+2b_{4}x+b_{6}$.

For positive integers $n$, we define the {\it division polynomials}
$\psi_{n}, \phi_{n} \in \bbZ \left[ a_{1}, a_{2}, a_{3},a_{4}, a_{6}\right] \left[ x, y \right]$ by
\begin{align*}
\psi_{1} &= 1, \\
\psi_{2} &= 2 y+a_{1}x+a_{3}, \\
\psi_{3} &= 3 x^{4}+b_{2}x^{3}+3 b_{4}x^{2}+3 b_{6}x+b_{8}, \\
\psi_{4} &= \psi_{2} \left( 2x^{6}+b_{2}x^{5}+5 b_{4}x^{4}+10 b_{6}x^{3}
            +10 b_{8}x^{2} + \left( b_{2}b_{8}-b_{4}b_{6} \right)x+b_{4}b_{8}-b_{6}^{2} \right),
\end{align*}
and then inductively by the formulas
\begin{align*}
\psi_{2m+1}       &= \psi_{m+2}\psi_{m}^{3}-\psi_{m-1}\psi_{m+1}^{3} &\text{for $m \geq 2$}, \\
\psi_{2}\psi_{2m} &= \psi_{m} \left( \psi_{m+2}\psi_{m-1}^{2}-\psi_{m-2}\psi_{m+1}^{2} \right) &\text{for $m \geq 3$} \nonumber.
\end{align*}

We also have
\begin{align}
\label{eq:phi1}
\phi_{1} &= x, \nonumber \\
\phi_{n} &= x \psi_{n}^{2}-\psi_{n-1}\psi_{n+1} &\text{for $n \geq 2$}.
\end{align}

We will sometimes need expressions for these polynomials that depend
only on $x$:
\begin{align}
\label{eq:phi2}
\phi_{2}(x)     &= x^{4}-b_{4}x^{2}-2b_{6}x-b_{8}, \\
\label{eq:psi2-sqrB}
\psi_{2}^{2}(x) &= 4x^{3}+b_{2}x^{2}+2b_{4}x+b_{6}.
\end{align}

In what follows, we will often use $P$ as the argument of these polynomials,
rather than $x(P)$ and $y(P)$.

For a finite extension, $L$, of $K$, we will use $R_{L}$, $\cM_{L}$ and $v_{L}$
for the ring of integers of $L$, its maximal ideal and the associated valuation,
respectively.

We put $c_{v}=\left| E(K)/E_{0}(K) \right|$,
where $E_{0}(K)=\left\{ P \in E(K): \widetilde{P} \in \widetilde{E}(k)_{\rm ns} \right\}$,
the set of points of $E(K)$ with non-singular reduction modulo $\pi$. This
quotient group is known as the component group of $E$ at $v$, so $c_{v}$ is
the order of the component group.

Tate's Algorithm \cite{Tate1} to compute the special fibre of a N\'{e}ron model
will play a crucial role in many parts of our work. We will use Silverman's
presentation of it in \cite[Chapter~IV, Section~9]{Silv3}.


\subsection{Simplifying $R_{n}(a,\ell)$}

\begin{definition}
\label{def:Rn}
To any pair $(a, \ell)$ of integers satisfying $\ell>0$, we associate an
integer sequence, $\left\{ R_{n}(a,\ell) \right\}_{n \geq 0}$, defined by
\[
R_{n}(a,\ell)
= \frac{n^{2} \widehat{a} \left( \ell -\widehat{a} \right)-\widehat{na} \left( \ell-\widehat{na} \right)}{2\ell},
\]
where $\widehat{x}$ denotes the least non-negative residue of $x$ modulo $\ell$.
\end{definition}

This sequence is identical to the sequence $R_{n}(a,\ell)$ defined in
Definition~8.1 of \cite{Stange}, but the expression here is simpler. We prove
that now.

\begin{lemma}
Let $a,\ell,n$ be non-negative integers with $\ell \geq 1$ and let $\hat{x}$
denote the least non-negative residue of $x$ modulo $\ell$. Then
\[
n^{2}\widehat{a} \left( \ell-\widehat{a} \right)
\equiv \widehat{na} \left( \ell-\widehat{na} \right) \bmod{2}\ell.
\]

As a consequence, $R_{n}(a,\ell)$ here is identical to $R_{n}(a,\ell)$ in
Definition~$8.1$ of \cite{Stange}.
\end{lemma}

\begin{proof}
Write $a=a_{1}\ell+a_{2}$ where $0 \leq a_{2}<\ell$ (i.e., $\widehat{a}=a_{2}$).
Then we can write $na=\left( na_{1}+a_{3} \right)\ell+na_{2}-a_{3}\ell$ with
$0 \leq na_{2}-a_{3}\ell<\ell$ for some integer $a_{3}$.

Thus
\begin{align*}
n^{2}\widehat{a} \left( \ell-\widehat{a} \right)
&= n^{2}a_{2} \left( \ell-a_{2} \right)=n^{2}a_{2}\ell-n^{2}a_{2}^{2}, \\
\widehat{na} \left( \ell-\widehat{na} \right)
&= \left( na_{2}-a_{3}\ell \right) \left( \ell-na_{2}+a_{3}\ell \right)
= a_{2}n\ell-a_{2}^{2}n^{2}+2a_{2}a_{3}n\ell-a_{3}^{2}\ell^{2}-a_{3}\ell^{2}.
\end{align*}

Subtracting these two expressions, we obtain
\[
n^{2}\widehat{a} \left( \ell-\widehat{a} \right) - \widehat{na} \left( \ell-\widehat{na} \right)
= \left( n^{2}-n \right)a_{2}\ell-2a_{2}a_{3}n\ell+\left( a_{3}^{2}+a_{3} \right)\ell^{2}.
\]

Since $x^{2} \pm x$ is even for any integer $x$, the congruence in the lemma holds.

The simpler expression for $R_{n}(a,\ell)$ is immediate.
\end{proof}


\subsection{Stange's results for $v\left( \psi_{n}(P) \right)$}
\label{sec:stange}

We next state some theorems from Stange's paper \cite{Stange}.
Following Definition~5.3 of \cite{Stange}, we let
\begin{equation}
\label{eq:param-conditions}
b \in p\bbZ^{>0} \cup \{ 1 \}, \hspace*{2.0mm} e \in \bbZ^{>0}, \hspace*{2.0mm}
h \in \bbZ^{\geq 0}, \hspace*{2.0mm} j \in \bbZ^{\geq 0}, \hspace*{2.0mm}
s \in \bbZ^{>0} \cup \{ \infty \}, \hspace*{2.0mm} w \in \bbZ^{\geq 0} \cup \{ \infty \},
\end{equation}
and for $n \in \bbZ$, put
\begin{equation}
\label{eq:s-defn}
S_{n}(p,b,e,h,s,w) = \left\{
	\begin{array}{ll}
		b^{j}s + \dfrac{b^{j}-1}{b-1}h + e \left( v_{p}(n) - j \right) + w & v_{p}(n) > j, \\
		b^{v_{p}(n)}s + \dfrac{b^{v_{p}(n)}-1}{b-1}h & v_{p}(n) \leq j,
	\end{array} \right.
\end{equation}
where $v_{p}$ is the valuation on $\bbQ$ associated to the rational prime, $p$
(recall our notation in Section~1).

\vspace*{1.0mm}

We shall use the following values for these quantities here.\\
--$b$ will be the smallest exponent of $T$ with the valuation $v$ of its coefficient less
than $v(p)$ in the expansion of multiplication-by-$p$, $[p]T$, in the formal
group of the elliptic curve $E$, or else $b=1$, if no such coefficient exists.

\noindent
--$e=v(p)$.

\noindent
--$h$ will be the valuation of the coefficient of $T^{b}$, or else $h=0$, if $b=1$.

\noindent
--$n_{P}$ will be the smallest positive integer such that
$\widetilde{\left[ n_{P} \right]P}=\widetilde{\cO}$, where $\widetilde{\cdot}$
is the reduction map from $E(K)$ to $\widetilde{E}(k)$ and $k$ is the
residue field, $k$.

\noindent
--$s=s_{P}=v \left( x \left( \left[ n_{P} \right] P \right) / y \left( \left[ n_{P} \right] P \right) \right)$.

\noindent
--$j=j_{P}=0$ if $b=1$; otherwise $j$ will be the smallest non-negative integer
such that
\[
e \leq b^{j}((b-1)s+h).
\]

\noindent
--$w=w_{P}=0$ unless $b>1$ and we have equality in the definition of $j$ above.
In this case (i.e., $b>1$ and we have equality in the definition of $j$), put
\begin{equation}
\label{eq:w-defn}
w = v \left( \frac{x \left( \left[ p^{j+1} n_{P} \right] P \right)}{y \left( \left[ p^{j+1} n_{P} \right] P \right)} \right)
- b v \left( \frac{x \left( \left[ p^{j} n_{P} \right] P \right)}{y \left( \left[ p^{j} n_{P} \right] P \right)} \right)
- h,
\end{equation}
which may be equal to $+\infty$.

To simplify our notation in what follows, we will often write $S_{n}(P)$ instead
of $S_{n}\left( p,b,v(p),h,s_{P},w_{P} \right)$.

\begin{remark}
(i) In the expression for $w$ in \eqref{eq:w-defn}, we have $b$ as the coefficient
of the second term. This corrects an error in the expression for $w$ in Lemma~5.1(iii)
of \cite{Stange}, where the exponent $p$ should be $b$.

\vspace*{1.0mm}

(ii) We have $S_{n} \in \bbZ^{>0} \cup \{ \infty \}$.

\vspace*{1.0mm}

(iii) In keeping with the conventional notation of $e$ for the ramification
index, we use $e$ here, where Stange has used $d$.

\vspace*{1.0mm}

(iv) When $b=1$, we have $h=0$ and $j=0$. In this case, we
use the convention $\left( b^{0}-1 \right)/(b-1)h=0$ to avoid the indeterminate
form $\left( b^{0}-1 \right)/(b-1)=0/0$ in the above expressions for $S_{n}(P)$.
\end{remark}

\begin{lemma}{\rm (\cite{Stange}, Theorem~6.1)}
\label{lem:non-sing}
Assume that $E$ is in minimal Weierstrass form and $P$ has non-singular reduction.
Then
\[
v \left( \psi_{n}(P) \right)
= \min\left\{ 0, \frac{v(x(P))}{2} \right\}n^{2}
+ \left\{
	\begin{array}{ll}
		S_{n/n_{P}}(P) & \text{if $n_{P} \mid n$,} \\
		0              & \text{if $n_{P} \nmid n$.}
	\end{array}
\right.
\]

Furthermore, $v(x(P))<0$ if and only if $n_{P}=1$.
\end{lemma}

\begin{remark}
From the definition of $n_{P}$, $v \left( x \left( [n]P \right) \right)<0$
if and only if $n_{P}|n$. Therefore if $n_{P} \nmid n$, then
$v\left( \phi_{n}(P) \right) \geq v\left( \psi_{n}^{2}(P) \right)$.
\end{remark}

\begin{lemma}{\rm (\cite{Stange}, Theorem~9.3)}
Suppose that $E$ is in minimal Weierstrass form with multiplicative reduction,
$P$ has singular reduction, and let $n_{P}$ be as above. Then
\[
v \left( \psi_{n}(P) \right) =
R_{n} \left( a_{P},\ell_{P} \right) + \left\{
\begin{array}{ll}
	S_{n/n_{P}} \left( p, p, v(p), 0, s_{P},w_{P} \right) & n_{P} \mid n, \\
	0 & n_{P} \nmid n,
\end{array} \right.
\]
where $\ell_{P}=-v(j(E))$ and $a_{P}$ is the component of the N\'{e}ron model
special fibre $(\cong \bbZ/\ell_{P}\bbZ)$ containing $P$.
\end{lemma}

\begin{lemma}
\label{lem:phi}
{\rm (i)} Suppose that $E$ is in minimal Weierstrass form, $P$ has non-singular
reduction, and $n_{P}$ is as above. If $n_{P} | n$ or $v \left( x \left( [n]P \right) \right)=0$,
then
\[
v\left( \phi_{n}(P) \right)=\min \left\{ 0, v(x(P)) \right\}n^{2}.
\]

{\rm(ii)} Suppose that $E$ is in minimal Weierstrass form with multiplicative
reduction and $P$ has singular reduction. If $n_{P} | n$, then
\[
v\left( \phi_{n}(P) \right)=2R_{n} \left( a_{P}, \ell_{P} \right),
\]
where $\ell_{P}=-v(j(E))$ and $a_{P}$ is the component of the N\'{e}ron model
special fibre $(\cong \bbZ/\ell_{P}\bbZ)$ containing $P$.
\end{lemma}

\begin{proof}
(i) Suppose that $n_{P} | n$. From Lemma~\ref{lem:non-sing},
\[
v \left( \psi_{n}(P) \right)
= \min\left\{ 0, \frac{v(x(P))}{2} \right\}n^{2} +
S_{n/n_{P}} (P).
\]

We also have
\[
v\left( \phi_{n}(P) \right) = v \left( \psi_{n}^{2}(P) \right)+v(x([n]P)).
\]

Since $n_{P} |n$, from the proof of Theorem~6.1 of \cite{Stange}, we have
$v \left( -x([n]P)/y([n]P) \right)=S_{n/n_{P}}(P)$. From the minimal Weierstrass
equation for $E$, we have $3v(x([n]P))=2v(y([n]P))$, since $v(x([n]P))<0$.
Hence $v(x([n]P))=-2v\left( -x([n]P)/y([n]P) \right)$ and so
\[
v\left( \phi_{n}(P) \right) = v \left( \psi_{n}^{2}(P) \right)
-2S_{n/n_{P}} (P).
\]

It follows from our expression above for $v \left( \psi_{n}(P) \right)$ that
$v\left( \phi_{n}(P) \right)=\min \left\{ 0, v(x(P)) \right\}n^{2}$.

If $v \left( x \left( [n]P \right) \right)=0$, then
\[
v\left( \phi_{n}(P) \right) = v \left( \psi_{n}^{2}(P) \right)+v(x([n]P))
v\left( \phi_{n}(P) \right) = v \left( \psi_{n}^{2}(P) \right).
\]

So if $n_{P} \nmid n$, then (i) follows immediately from Lemma~\ref{lem:non-sing}.

(ii) Suppose that $n_{P} | n$. From the proof of Lemma~11.4 in \cite{Stange},
\[
v \left( \psi_{n}(P) \right)=R_{n} \left( a_{P}, \ell_{P} \right) + v \left( x([n]P)/y([n]P) \right).
\]
 
Therefore,
\begin{align*}
v \left( \phi_{n}(P) \right)
&= v \left( \psi_{n}^{2}(P) \right) + v(x([n]P)) \\
&= 2R_{n} \left( a_{P}, \ell_{P} \right) + 2v \left( x([n]P)/y([n]P) \right)+v(x([n]P)).
\end{align*}

Since $3v(x([n]P))=2v(y([n]P))$, we have
\[
v \left( x([n]P)/y([n]P) \right) = v(x([n]P))-v(y([n]P))=-v(x([n]P))/2.
\]

Hence $v \left( \phi_{n}(P) \right) = 2R_{n} \left( a_{P}, \ell_{P} \right)$.
\end{proof}


\subsection{Non-integral $x(P)$}
\label{subsect:non-int}

By Proposition~2.2(iii) of \cite{Stange}, if $E$ is given by a $v$-integral Weierstrass
equation, where $v$ is a nonarchimedean valuation, and $v(x(P)),v(y(P)) < 0$,
then $v \left( \phi_{n}(P) \right) = n^{2}v(x(P))$. From the minimal Weierstrass
equation for $E$, we have $3v(x(P)) = 2v(y(P))$, so it suffices that $v(x(P))<0$.

The degree of $\psi_{n}^{2}(P)$ as a polynomial in $x(P)$ is $n^{2}-1$ (see
Proposition~2.2(ii) of \cite{Stange}, for example). So if $E$ is given by a $v$-integral Weierstrass
equation, where $v$ is a nonarchimedean valuation, and $v(x(P)) < 0$, then
$v \left( \psi_{n}^{2}(P) \right) \geq \left( n^{2}-1 \right) v(x(P))$.

Thus $k_{v,n}(P)= n^{2}v(x(P))$ in this case. Also notice that if $v(x(P))<0$,
then $\widetilde{P}=\widetilde{\cO}$ and so $P$ is non-singular modulo $\pi$.
This establishes Theorem~\ref{thm:main} in this case.


\subsection{Points with singular reduction}

The following lemma will allow us to use the Weierstrass equations obtained in
the course of Tate's Algorithm to simplify our work. Throughout this section,
we will let $P'=(x',y')$ be the image of $P=(x,y)$ under a change of variables
of the form $x=u^{2}x'+r$ and $y=u^{3}y'+u^{2}sx'+t$ with $r,s,t,u \in R$.

\begin{lemma}
\label{lem:sing-pts}
{\rm (i)} Let $P=(x,y) \in K^{2}$ satisfy the equation
\[
f(x,y)=y^{2}+a_{1}xy+a_{3}y-x^{3}-a_{2}x^{2}-a_{4}x-a_{6}=0,
\]
with $a_{1}, a_{2}, a_{3}, a_{4}, a_{6} \in R$ and $v\left( a_{3} \right)>0$,
$v\left( a_{4} \right)>0$ and $v\left( a_{6} \right)>0$.

$P$ has singular reduction if and only if $v(x)>0$ and $v(y)>0$.

{\rm (ii)} Let $n_{P}$ be as defined in Subsection~$\ref{sec:stange}$. If
$v(u)=0$, then $n_{P'}=n_{P}$.

{\rm (iii)} Suppose that $E$ has multiplicative reduction and $P \in E(K)$
has singular reduction. Let $a_{P}$ be the component of the cyclic group
$E(K)/E_{0}(K)$ that contains $P$. If $v(u)=0$, then $a_{P'}=a_{P}$.

{\rm (iv)} If $v(u)=0$, then $k_{v,n}(P)=k_{v,n}(P')$. As a consequence, if $E$ is
in minimal Weierstrass form, then the changes of variable in Tate's Algorithm
leave $k_{v,n}(P)$ unchanged.
\end{lemma}

\begin{proof}
(i)
\begin{align}
\label{eq:singPt-2}
\text{$P$ has singular reduction}
&\Longleftrightarrow v \left( \frac{\partial f}{\partial x}(P) \right)>0
\text{ and } v \left( \frac{\partial f}{\partial y}(P) \right)>0 \notag \\
&\Longleftrightarrow
\begin{cases}
v \left( a_{1}y-3x^{2}-2a_{2}x \right)>0 \text{ and} \\
v \left( 2y+a_{1}x \right)>0.
\end{cases}
\end{align}
The last logical equivalence here holds using the expressions for $\partial f/\partial x(P)$
and $\partial f/\partial y(P)$ since we assume that $v\left( a_{3} \right)>0$
and $v\left( a_{4} \right)>0$.

If $v(x)>0$ and $v(y)>0$, then it is immediate that $P$ has singular reduction.
So we may assume that $P$ has singular reduction.

Subtracting $2$ times the first expression in \eqref{eq:singPt-2} from $a_{1}$
times the second one, we can eliminate $y$ and obtain
\[
v \left( x \left( 6x+a_{1}^{2}+4a_{2} \right) \right) >0.
\]

Similarly, using the second expression in \eqref{eq:singPt-2} to eliminate $y$
from $f(x,y)$, we have
\[
v \left( x^{2} \left( 4x+a_{1}^{2}+4a_{2} \right) \right) > 0,
\]
since $v \left( a_{3} \right)>0$ and $v \left( a_{4} \right)>0$.

If $v(x)>0$, then from $f(x,y)=0$, $v \left( a_{3} \right)>0$ and
$v \left( a_{6} \right)>0$, we have $v(y)>0$.

Assume that $v(x) \leq 0$. Then
\[
v \left( 4x+a_{1}^{2}+4a_{2} \right)>0
\quad \text{and} \quad v \left( 6x+a_{1}^{2}+4a_{2} \right)>0.
\]

Subtracting these two expressions, we obtain $v(2x)>0$, which contradicts $v(x) \leq 0$
when $v(2)=0$.

If $v(2)>0$, then from \eqref{eq:singPt-2}, we have
$v \left( a_{1}y-3x^{2} \right)>0$ and $v \left( a_{1}x \right)>0$. Subtracting
$y$ times the second expression from $x$ times the first expression, we obtain
$v \left( 3x^{2} \right)>0$, so $v(x)>0$. As above, $v(x)>0$ also
implies that $v(y)>0$.

(ii) We know that $\widetilde{P}=\widetilde{\cO}$ if and only if $v(x)<0$.

Suppose that $v(x')<0$. We have $v(x)=v \left( u^{2}x'+r \right) \geq \min \left\{ v(x'), v(r) \right\}$
with equality if $v(x') \neq v(r)$. Since $r \in R$, we have $v(r) \geq 0$. Therefore
$v(x)=v(x')<0$ and $n_{P}|n_{P'}$. The same argument also shows that $n_{P'}|n_{P}$,
so $n_{P}=n_{P'}$.

(iii) This follows from the expression for $a_{P}$ in Lemma~5.1 of \cite{Silv4}
(denoted there as $n$), using the facts that $\Delta=u^{12}\Delta'$,
$2y+a_{1}x+a_{3}=u^{3} \left( 2y'+a_{1}'x'+a_{3}' \right)$ and $v(u)=0$.

(iv) By Proposition~2.2(iv) of \cite{Stange}, our change of variables gives
\begin{equation}
\label{eq:psiN}
\psi_{n}(P) = u^{n^{2}-1}\psi_{n}\left( P' \right).
\end{equation}

From this and \eqref{eq:phi1}, we obtain
\begin{align}
\label{eq:phiN}
\phi_{n}(P) &= x \psi_{n}^{2}(P)-\psi_{n-1}(P)\psi_{n+1}(P) \nonumber \\
&= \left( u^{2}x\left( P' \right)+r \right) u^{2n^{2}-2} \psi_{n}^{2}\left( P' \right)
   - u^{2n^{2}}\psi_{n-1}\left( P' \right)\psi_{n+1}\left( P' \right) \nonumber \\
&= u^{2n^{2}-2} \left( u^{2}\phi_{n}\left( P' \right)+r\psi_{n}^{2}\left( P' \right) \right).
\end{align}

From \eqref{eq:psiN} and $v(u)=0$, we have $v \left( \psi_{n}^{2}(P') \right) =v \left( \psi_{n}^{2}(P) \right)$.

If $v \left( \phi_{n}(P) \right)
< v \left( \psi_{n}^{2}(P) \right)=v \left( \psi_{n}^{2}(P') \right)$, then $v \left( \phi_{n}(P') \right)
=v \left( \phi_{n}(P) \right)$ from \eqref{eq:phiN}. So $k_{v,n}(P')=v \left( \phi_{n}(P') \right)
=k_{v,n}(P)$.

Suppose that $v \left( \phi_{n}(P) \right) \geq v \left( \psi_{n}^{2}(P) \right)
=v \left( \psi_{n}^{2}(P') \right)$.
If $v \left( \phi_{n}\left( P' \right) \right) \geq v \left( r\psi_{n}^{2}\left( P' \right) \right)$,
then $k_{v,n}(P')=v \left( \psi_{n}^{2}\left( P' \right) \right)$, as desired.
If $v \left( \phi_{n}\left( P' \right) \right)<\left( r\psi_{n}^{2}\left( P' \right) \right)$,
then $v \left( \phi_{n}(P) \right)=v \left( \phi_{n}\left( P' \right) \right)$,
from \eqref{eq:phiN}. So $v \left( \phi_{n}\left( P' \right) \right)
\geq v \left( \psi_{n}^{2}(P') \right)$, as required.

The statement regarding the changes of variables in Tate's Algorithm now
follows because it is only in Step~11 (i.e., when we do not start with a
minimal model) that we perform a change of variables with $v(u) \neq 0$.
\end{proof}

\subsection{Additive reduction}

In order to use Stange's results above, we will often need to work in a finite
extension of $K$. Here we present results on how the valuations behave when we
work in such extensions. In addition, in Lemmas~\ref{lem:ImStar-vals} and
\ref{lem:pot-good-red-vals} we provide some results we will need for the valuation
of $\phi_{n}(P)$ and $\psi_{n}(P)$ in order to finish the proof of Theorem~\ref{thm:main}.

\begin{lemma}
\label{lem:val-extensions}
Let $E/K$ be an elliptic curve having additive reduction and $P \in E(K)$ having
singular reduction.

Let $x=u^{2}x'+r$ and $y=u^{3}y'+u^{2}sx'+t$ be a change of variables from $E$
to an elliptic curve $E'$ in minimal Weierstrass form with $r,s,t,u \in R_{L}$,
where $L$ is a finite extension of $K$. Write $P'$ for the image of $P$ under
this change of variables. Put
\[
T(E) = \left\{
\begin{array}{ll}
	v_{L} \left( \Delta_{E} \right)/6 & \text{if $E'$ has good reduction}, \\
	v_{L} \left( c_{4}(E) \right)/2   & \text{if $E'$ has multiplicative reduction}.
\end{array}
\right.
\]

{\rm (i)} We have $v_{L}(u)=T(E)/2$ and
\begin{align}
\label{eq:psi-extension}
v_{L} \left( \psi_{n}(P) \right)
&= \left( n^{2}-1 \right)T(E)/2+v_{L} \left( \psi_{n}\left( P' \right) \right), \\
\label{eq:phi-extension}
v_{L} \left( \phi_{n}(P) \right)
&= \left( n^{2}-1 \right)T(E) + v_{L} \left( u^{2} \phi_{n}\left( P' \right)+r \psi_{n}^{2}\left( P' \right) \right).
\end{align}

{\rm (ii)} If $[n]P$ has non-singular reduction, $v \left( a_{3} \right)>0$, $v \left( a_{4} \right)>0$
and $v \left( a_{6} \right)>0$, then
\begin{equation}
\label{eq:lemma 13}
v_{L} \left( \phi_{n}(P) \right)
= \left( n^{2}-1 \right) T(E)+ v_{L} \left( u^{2} \phi_{n}\left( P' \right) \right).
\end{equation}
\end{lemma}

\begin{proof}
(i) Equation~\eqref{eq:psi-extension} follows from \cite{Stange}. We use Theorem~7.1
when $E$ has potential good reduction and Theorem~9.3 when $E$ has potential
multiplicative reduction.

For Equation~\eqref{eq:phi-extension}, we use \eqref{eq:phiN} in the proof of
Lemma~\ref{lem:sing-pts}(iv):
\[
\phi_{n}(P) = u^{2n^{2}-2} \left( u^{2}\phi_{n}\left( P' \right)+r\psi_{n}^{2}\left( P' \right) \right).
\]

Suppose that $E'$ has good reduction. Since $\Delta_{E'}=u^{-12}\Delta_{E}$, we
have $v_{L}\left( \Delta_{E'} \right)=v_{L}\left( \Delta_{E} \right)-12v_{L}(u)$.
It must equal $0$, since $E'$ has good reduction. Therefore,
$v_{L}(u)=v_{L}\left( \Delta_{E} \right)/12$.

Now suppose that $E'$ has multiplicative reduction. Since $c_{4}\left(E'\right)=u^{-4}c_{4}(E)$,
we have $v_{L}\left( c_{4}\left(E'\right) \right)=v_{L}\left( c_{4}(E) \right)-4v_{L}(u)$.
By \cite[Proposition~VII.5.1(b)]{Silv2}, $v_{L}\left( c_{4}\left(E'\right) \right)=0$,
since $E'$ has multiplicative reduction. Therefore, $v_{L}(u)=v_{L}\left( c_{4}(E) \right)/4$.

(ii) Putting
$V_{P', n}=v_{L} \left( u^{2} \phi_{n}\left( P' \right)+r \psi_{n}^{2} \left(P'\right) \right)$,
we will show that $V_{P',n}=v_{L}\left(u^{2}\phi_{n}\left(P'\right)\right)$.

Since $[n]P$ has non-singular reduction, by our assumptions on the $a_{i}$'s and
Lemma~\ref{lem:sing-pts}(i), we have
$v_{L}\left( \phi_{n}(P) \right) \leq v_{L}\left( \psi_{n}^{2}(P) \right)$.

From \eqref{eq:phi-extension}, we have
$V_{P',n}=v_{L} \left( \phi_{n}(P) \right)-\left( n^{2}-1 \right)T(E)$.
Applying $v_{L}\left( \phi_{n}(P) \right) \leq v_{L}\left( \psi_{n}^{2}(P) \right)$
and then \eqref{eq:psi-extension}, we obtain
\begin{equation}
\label{eq:VP'}
V_{P',n} \leq v_{L} \left( \psi_{n}^{2} \left(P'\right) \right).
\end{equation}

We now consider the case when $v_{L}(r)>0$.
Combining $V_{P',n} \geq \min \left( v_{L} \left( u^{2} \phi_{n}\left( P' \right) \right),
v_{L} \left( r \psi_{n}^{2} \left(P'\right) \right) \right)$ with \eqref{eq:VP'},
since $v_{L}(r)>0$, it follows that $V_{P',n}=v_{L}\left( u^{2} \phi_{n}\left( P' \right)\right)$.

Next we consider the case when $v_{L}(r)=0$. Since $P$ has singular reduction,
by Lemma~\ref{lem:sing-pts}(i), we have $v_{L}(x(P))>0$. Since $v_{L}(r)=0$,
we have $v_{L}(x(P)-r)=0$. Therefore
\[
v_{L} \left( x \left(P'\right) \right) = v_{L} \left( u^{-2}(x(P)-r) \right)
=-v_{L} \left( u^{2} \right)<0,
\]
the last inequality holding by \cite[Proposition~VII.5.1(c)]{Silv2}, our definition
of $T(E)$ and the relationship between $T(E)$ and $u$.

Hence $v_{L} \left( x\left( [n]P' \right) \right) \leq v_{L} \left( x\left( P' \right) \right)<0$.
The first inequality here follows from the our arguments in Subsection~\ref{subsect:non-int},
which show that $v \left( \phi_{n}(P') \right)=n^{2}v \left( x(P') \right)$ and
$v \left( \psi_{n}^{2}(P') \right) \geq \left( n^{2}-1 \right)v \left( x(P') \right)$.
Hence $v_{L}\left( \phi_{n} \left(P'\right) \right) - v_{L}\left( \psi_{n}^{2}\left(P'\right)\right)
\leq -v_{L} \left( u^{2} \right)$, which implies that $v_{L} \left( u^{2} \phi_{n}\left( P' \right)\right)
\leq v_{L}\left( \psi_{n}^{2}\left(P'\right) \right) = v_{L} \left( r\psi_{n}^{2} \left( P' \right) \right)$.
If $v_{L} \left( u^{2} \phi_{n}\left( P' \right)\right) < v_{L} \left( r\psi_{n}^{2} \left( P' \right) \right)$
then $V_{P',n}=v_{L} \left( u^{2} \phi_{n}\left( P' \right) \right)$. If
$v_{L} \left( u^{2} \phi_{n}\left( P' \right)\right) = v_{L}\left( r \psi_{n}^{2} \left( P' \right) \right)$,
then $V_{P',n} \geq v_{L} \left( r\psi_{n}^{2} \left( P' \right) \right)$. Combining
this with $V_{P',n} \leq v_{L} \left( r\psi_{n}^{2} \left(P'\right) \right)$
(from \eqref{eq:VP'} and $v_{L}(r)=0$), we obtain $V_{P',n}=v_{L} \left( r\psi_{n}^{2} \left(P'\right) \right)
=v_{L} \left( u^{2} \phi_{n}\left( P' \right) \right)$.
\end{proof}

\begin{lemma}
\label{lem:ImStar-vals}
Let $E/K$ be an elliptic curve having Kodaira type $I_{m}^{*}$ with $m \geq 1$, and $P \in E(K)$ having singular reduction.
By Tate's algorithm, we may assume that the Weierstrass equation for $E$ has
$v \left( a_{1} \right) \geq 1$, $v \left( a_{2} \right)=1$, $v \left( a_{3} \right) \geq \lfloor m/2 \rfloor+2$,
$v \left( a_{4} \right) \geq \lfloor (m-1)/2 \rfloor+3$ and $v \left( a_{6} \right) \geq m+3$.

{\rm (i)} Let $m$ be odd. If $v\left(x(P)\right)=1$, then $[2]P$ has non-singular
reduction. Here $v\left( \phi_{2}(P) \right)=v\left( \psi_{3}(P) \right)=4$
and $v \left( \psi_{2}^{2}(P) \right) \geq 4$.
$1<v \left(x(P) \right) < (m+3)/2$ is not possible.
If $v \left(x(P) \right) \geq (m+3)/2$, then $[2]P$ has singular reduction.
Here $v\left( \phi_{2}(P) \right)=v\left( \psi_{3}(P) \right)=m+4$
and $v \left( \psi_{2}^{2}(P) \right)=m+3$.

{\rm (ii)} If $m$ is even, then $[2]P$ has non-singular reduction. If
$v\left(x(P)\right)=1$, then $v\left( \phi_{2}(P) \right)=v\left( \psi_{3}(P) \right)=4$
and $v \left( \psi_{2}^{2}(P) \right) \geq 4$.
$1<v \left(x(P) \right) < (m+2)/2$ is not possible.
If $v \left(x(P) \right) \geq (m+2)/2$, then
$v\left(\phi_{2}(P) \right)=v\left( \psi_{3}(P) \right)=m+4$
and $v \left( \psi_{2}^{2}(P) \right)=m+4$.
\end{lemma}

\begin{proof}
We write $x$ in place of $x(P)$ for convenience in what follows. Since
$v \left( a_{3} \right)>0$, $v \left( a_{4} \right)>0$, $v \left( a_{6} \right)>0$
and $P$ has singular reduction, from Lemma~\ref{lem:sing-pts}(i), we have $v(x)>0$.

(i) Since $m$ is odd, write $m=2k-1$ for some positive integer $k$. From the
inequalities in the lemma, we have $v \left( a_{1} \right) \geq 1$,
$v \left( a_{2} \right)=1$, $v \left( a_{3} \right) \geq k+1$,
$v \left( a_{4} \right) \geq k+2$ and $v \left( a_{6} \right) \geq 2k+2$. In
addition, from Proposition~1(a) in Section~III of \cite{Papa} (note that $n$
there equals our $k+1$), we have
$v \left( b_{6} \right)=2k+2$ and $v \left( b_{8} \right)=2k+3$. As a consequence,
for $p \geq 3$, we also have $v \left( b_{2} \right)=1$ and $v \left( b_{4} \right) \geq k+2$.
For $p=2$, we have $v \left( b_{2} \right) \geq 2$ and $v \left( b_{4} \right) \geq k+2$.
Furthermore, from the proof of Proposition~1(a) in Section~III of \cite{Papa},
we find that $v \left( a_{3} \right)=k+1$ when $m$ is odd and $p=2$.

(i-a) Let $p \geq 3$.
For use with the expression for $\phi_{2}(P)$ in \eqref{eq:phi2}, we have
$v \left( x^{4} \right)=4v(x)$, $v \left( b_{4}x^{2} \right) \geq 2v(x)+k+2$,
$v \left( 2b_{6}x \right)=v(x)+2k+2$ and $v \left( b_{8} \right)=2k+3$. Similarly,
for use in the expression for $\psi_{2}^{2}(P)$ in \eqref{eq:psi2-sqrB}, we
have $v \left( 4x^{3} \right)=3v(x)$, $v \left( b_{2}x^{2} \right)=2v(x)+1$,
$v \left( 2b_{4}x \right) \geq v(x)+k+2$ and $v \left( b_{6} \right)=2k+2$.

If $v(x)=1$, then these inequalities imply that $v\left(\phi_{2}(P)\right)=v \left( x^{4} \right)=4$
and $v\left( \psi_{2}^{2}(P) \right) \geq 4$ (the expression for $\psi_{2}^{2}(P)$
in \eqref{eq:psi2-sqrB} shows that $v \left( \psi_{2}^{2}(P) \right) \geq v \left( 4x^{3} \right)=3$,
but $v \left( \psi_{2}^{2}(P) \right)$ is even),
so $[2]P$ has non-singular reduction.

If $2 \leq v(x) \leq k$, then from the inequalities above, we have
$v \left( 2b_{4}x \right) \geq v(x)+k+2 \geq 2v(x)+2$ and
$v \left( b_{6} \right)=2k+2 \geq 2v(x)+2$. Hence
$v \left( \psi_{2}^{2}(P) \right)=v \left( b_{2}x^{2} \right)=2v(x)+1$, which
is not possible.

If $v(x) \geq k+1=(m+3)/2$, then $v \left( \phi_{2}(P) \right)=v \left( b_{8} \right)=2k+3$
and $v \left( \psi_{2}^{2}(P) \right)=v \left( b_{6} \right)=2k+2$.
Thus $v \left( x\left([2]P \right) \right)=1$, so $[2]P$ is singular.

(i-b) Let $p=2$.
For use with the expression for $\phi_{2}(P)$ in \eqref{eq:phi2}, we have
$v \left( x^{4} \right)=4v(x)$, $v \left( b_{4}x^{2} \right) \geq 2v(x)+k+2$,
$v \left( 2b_{6}x \right) \geq v(x)+2k+3$ and $v \left( b_{8} \right)=2k+3$. For the
expression for $\psi_{2}^{2}(P)$ in \eqref{eq:psi2-sqrB}, we have
$v \left( 4x^{3} \right) \geq 3v(x)+2$, $v \left( b_{2}x^{2} \right) \geq 2v(x)+2$,
$v \left( 2b_{4}x \right) \geq v(x)+k+3$ and $v \left( b_{6} \right)=2k+2$.

We consider separately three cases according to the value of $v(x)$.

If $v(x)=1$, then we have $v\left( \phi_{2}(P) \right)=v \left( x^{4} \right)=4$
and $v\left( \psi_{2}^{2}(P) \right) \geq v \left( b_{2}x^{2} \right) \geq 4$,
so $[2]P$ has non-singular reduction.

\vspace*{1.0mm}

Next assume that $2 \leq v(x) \leq k$.
Using our inequalities for the
$v \left( a_{i} \right)$'s, we see that
$v \left( x^{3}+a_{2}x^{2}+a_{4}x+a_{6} \right) = v \left( a_{2}x^{2} \right)
=2v(x)+1$, since $v\left( a_{4}x \right)\geq k+2+v(x)\geq 2v(x)+2$ and
$v\left( a_{6} \right)\geq 2k+2 \geq 2v(x)+2$.
We also have $v \left( y^{2} \right)=2v(y)$, $v \left( a_{1}xy \right) \geq 1+v(x)+v(y)$,
$v \left( a_{3}y \right)=k+1+v(y) \geq v \left( a_{1}xy \right)$.

If $v(y) \leq v(x)$, then $v \left( y^{2}+a_{1}xy+a_{3}y \right)=v \left( y^{2} \right) \leq 2v(x)$.
If $v(y)>v(x)$, then $v \left( y^{2}+a_{1}xy+a_{3}y \right)>2v(x)+1$. In both
cases, the valuation of the left-hand side of the Weierstrass equation differs
from the valuation of the right-hand side, so $2 \leq v(x) \leq k$ is not possible.

\vspace*{1.0mm}

Lastly, assume that $v(x) \geq k+1$ with $k \geq 1$. Then $v\left(\phi_{2}(P)\right)
=v\left( b_{8} \right)=2k+3$ and $v\left( \psi_{2}^{2}(P) \right)=v\left( b_{6} \right)
=2k+2$, so $v\left( x\left( [2]P \right) \right)=1$. Therefore $[2]P$ has singular reduction.

\vspace*{1.0mm}

(ii) Since $m$ is even, write $m=2k$ for some positive integer $k$. From the
inequalities in the lemma, we have $v \left( a_{1} \right) \geq 1$,
$v \left( a_{2} \right)=1$, $v \left( a_{3} \right) \geq k+2$,
$v \left( a_{4} \right) \geq k+2$ and $v \left( a_{6} \right) \geq 2k+3$.
In addition, from Proposition~1(b) in Section~III of \cite{Papa}, we have $v \left( b_{8} \right)=2k+4$.
So, for $p \geq 3$, we have $v \left( b_{2} \right)=1$, $v \left( b_{4} \right) \geq k+2$
and $v \left( b_{6} \right) \geq 2k+3$. For $p=2$, we have $v \left( b_{2} \right) \geq 2$,
$v \left( b_{4} \right) \geq k+3$ and $v \left( b_{6} \right) \geq 2k+4$.
Furthermore, from the proof of Proposition~1(b) in Section~III of \cite{Papa},
we find that $v \left( a_{4} \right)=k+2$ when $m$ is even and $p=2$.

\vspace*{1.0mm}

(ii-a) Let $p \geq 3$.
For use with the expression for $\phi_{2}(P)$ in \eqref{eq:phi2}, we have
$v \left( x^{4} \right)=4v(x)$, $v \left( b_{4}x^{2} \right) \geq 2v(x)+k+2$,
$v \left( 2b_{6}x \right) \geq v(x)+2k+3$ and $v \left( b_{8} \right)=2k+4$.
For use with the expression for $\psi_{2}^{2}(P)$ in \eqref{eq:psi2-sqrB}, we
have $v\left( 4x^{3} \right) = 3v(x)$, $v \left( b_{2}x^{2} \right)=2v(x)+1$,
$v \left( 2b_{4}x \right) \geq v(x)+k+2$, and $v \left( b_{6} \right) \geq 2k+3$.

If $v(x)=1$, then $v \left( \phi_{2}(P) \right)=v \left( x^{4} \right)=4$ and
$v \left( \psi_{2}^{2}(P) \right) \geq 4$, so $[2]P$ is non-singular. If
$2 \leq v(x) \leq k$, then $v\left( \psi_{2}^{2}(P) \right)=v \left( b_{2}x^{2} \right)
=2v(x)+1$, which is not possible. If $v(x) \geq k+1$, then
$v \left( \psi_{2}^{2} (P)\right) \geq 2k+4$ and
$v\left( \phi_{2}(P) \right)=v \left( b_{8} \right)=2k+4$.
Therefore $[2]P$ has non-singular reduction.

\vspace*{1.0mm}

(ii-b) Let $p=2$.
From the expression for $\phi_{2}(P)$ in \eqref{eq:phi2}, we have $v \left( x^{4} \right)=4v(x)$,
$v \left( b_{4}x^{2} \right) \geq 2v(x)+k+3$, $v \left( 2b_{6}x \right) \geq v(x)+2k+5$
and $v \left( b_{8} \right)=2k+4$. From the expression for $\psi_{2}^{2}(P)$ in
\eqref{eq:psi2-sqrB}, we
have $v \left( 4x^{3} \right) \geq 3v(x)+2$, $v \left( b_{2}x^{2} \right) \geq 2v(x)+2$,
$v \left( 2b_{4}x \right) \geq v(x)+k+4$ and $v \left( b_{6} \right) \geq 2k+4$.

If $v(x)=1$, then $v \left( \phi_{2}(P) \right)=v \left( x^{4} \right)=4$ and
$v \left( \psi_{2}^{2}(P) \right) \geq 4$, so $[2]P$ has non-singular reduction.

Assume that $2 \leq v(x) \leq k$. We proceed in the same way as for such $v(x)$
in case~(i-b).
%

Assume that $v(x) \geq k+1$. We have $v \left( 4x^{3} \right) \geq 3k+5$,
$v \left( b_{2}x^{2} \right) \geq 2k+4$, $v \left( 2b_{4}x \right) \geq 2k+5$
and $v \left( b_{6} \right) \geq 2k+4$. Therefore $v \left( \psi_{2}^{2}(P) \right) \geq 2k+4$.
We also have $v \left( x^{4} \right) \geq 4k+4$, $v \left( b_{4}x^{2} \right) \geq 3k+5$,
$v \left( 2b_{6}x \right) \geq 3k+6$ and $v \left( b_{8} \right) =2k+4$. Therefore
$v \left( \phi_{2}(P) \right)=2k+4$ and hence $[2]P$ has non-singular reduction.
\end{proof}

\begin{lemma}
\label{lem:pot-good-red-vals}
Let $E/K$ be an elliptic curve. Suppose $P \in E(K)$ has singular reduction.
Upon applying Tate's algorithm, we have the following relationships.

{\rm (i)} If $m_{P}=2$, then $v \left( \phi_{2}(P) \right) =v \left( \psi_{3}(P) \right)$.

{\rm (ii)} If $m_{P}=3$ and $E/K$ has additive reduction,
then $v \left( \phi_{3}(P) \right)=3v \left(\psi_{2}^{2}(P) \right)$.
\end{lemma}

\begin{proof}
All assertions below for the values of the $v \left( a_{i} \right)$'s and
$v \left( b_{i} \right)$'s can be found in Silverman's presentation of Tate's
Algorithm in \cite[Chapter~IV, Section~9]{Silv3}. We will also use, sometimes
implicitly, the fact that $v(x(P))>0$. This follows from Lemma~\ref{lem:sing-pts}(i).

(i) Since $m_{P}=2$, it follows that $v \left( \phi_{2}(P) \right) \leq v \left( \psi_{2}^{2}(P) \right)$.
From Lemma~\ref{lem:sing-pts}(i), we have $v(x(P))>0$. So from the formula
$\phi_{2}=x \psi_{2}^{2}-\psi_{1}\psi_{3}$ (which follows from
\eqref{eq:phi1} with $n=2$) we obtain
$v \left(\phi_{2} (P) \right)=v \left( \psi_{1}(P) \psi_{3}(P) \right)=v \left( \psi_{3}(P) \right)$.

(ii) We again use \eqref{eq:phi1}, here with $n=3$, so we have
$\phi_{3}=x \psi_{3}^{2}-\psi_{2}\psi_{4}$.
Since $m_{P}=3$, it follows that
$v \left(\phi_{3}(P) \right) \leq v \left(\psi_{3}^{2}(P) \right)$
by Lemma~\ref{lem:sing-pts}(i).
Thus $v \left(\phi_{3}(P) \right)=v \left( \psi_{2}(P)\psi_{4}(P) \right)$, so
the result would follow if we showed that
$v \left( \psi_{4}(P) \right)=5v \left( \psi_{2}(P) \right)$.

Since $m_{P}=3$, it follows that $3|c_{v}$, so since $E/K$ has additive reduction,
either $E/K$ has Kodaira type $IV$ or Kodaira type $IV^{*}$.
We consider these two cases separately.

If $E/K$ has Kodaira type $IV$, then we are in Step~5 of Tate's Algorithm.
We have $v \left( b_{2} \right) \geq 1$, $v \left( b_{4} \right) \geq 2$,
$v \left( b_{6} \right) = 2$ and $v \left( b_{8} \right) \geq 3$. By using the
expression for $\psi_{2}^{2}(P)$ in \eqref{eq:psi2-sqrB}, we obtain $v \left(\psi_{2}^{2}(P) \right)=v \left( b_{6} \right)=2$.
From the expression for $\psi_{4}(P)$ in
Subsection~\ref{subsect:notation}, we have $v \left(\psi_{4}(P) \right)=
v \left(\psi_{2}(P) \right)+v \left( b_{6}^{2} \right)=5=5v \left(\psi_{2}(P) \right)$,
as desired.

If $E/K$ has Kodaira type $IV^{*}$, then we are in Step~8 of Tate's Algorithm, where $v \left( a_{1} \right) \geq 1$,
$v \left( a_{2} \right) \geq 2$, $v \left( a_{3} \right) \geq 2$,
$v \left( a_{4} \right) \geq 3$ and $v \left( a_{6} \right) \geq 4$. In addition,
$v \left( a_{3}^{2}/\pi^{4}+4a_{6}/\pi^{4} \right)
=v \left( a_{3}^{2}+4a_{6} \right)-4=0$.
If $p=2$, then $v \left( b_{2} \right) \geq 2$, $v \left( b_{4} \right) \geq 3$,
$v \left( b_{6} \right) = 4$ (from $v \left( a_{3}^{2}+4a_{6} \right)-4=0$ above)
and $v \left( b_{8} \right) \geq 6$.
If $p \geq 3$, then $v \left( b_{2} \right) \geq 2$, $v \left( b_{4} \right) \geq 3$,
$v \left( b_{6} \right) = 4$ and $v \left( b_{8} \right) \geq 4$.
Proceeding in the same way as for Kodaira type $IV$, we have
$v \left(\psi_{2}^{2}(P) \right)=v \left( b_{6} \right)=4$ and
$v \left(\psi_{4}(P) \right)=v \left(\psi_{2}(P) \right)+v \left( b_{6}^{2} \right)=10
=5v \left(\psi_{2}(P) \right)$.
\end{proof}


\section{Proof of Theorem~\ref{thm:main}}


\subsection{Multiplicative reduction}

Suppose that $E$ has multiplicative reduction and that $P$ is singular modulo $\pi$.

As in the proof of Lemma~11.3 of \cite{Stange}, we have
\[
k_{v,n}(P) = 2 v \left( \psi_{n}(P) \right) - 2 v \left( D_{n} \right)
= 2R_{n} \left( a_{P,v}, m \right),
\]
where $a_{P,v}$ is the component of the N\'{e}ron model special fibre
($\cong \bbZ/ m \bbZ$) containing $P$.

For non-split multiplicative reduction, we still have $k_{v,n}(P)
=2R_{n} \left( a_{P,v}, m \right)$ with $m$ even and $a_{P,v}=m/2$.

From Proposition~5(c) in \cite{CPS}, we see that the coefficient of $n^{2}$ in
$2R_{n} \left( a_{P,v}, m \right)$ is as stated in our theorem and in the
entry for $I_{m}$ in Table~\ref{table:eps}. The remaining term matches
$\epsilon_{v,n}(P)$ in the entry for $I_{m}$ in Table~\ref{table:eps} by a
simple elementary manipulation.


\subsection{Additive reduction}

Throughout this section, by Tate's algorithm, we may assume that the Weierstrass
equation for $E$ has $v \left( a_{3} \right)>0$, $v \left( a_{4} \right)>0$ and
$v \left( a_{6} \right)>0$.

By Proposition~5.5 in Chapter~VII of \cite{Silv2}, $E$ has potential good reduction
if and only if its $j$-invariant is integral. By Tableaux $I$, $II$ and $IV$ of
\cite{Papa}, we compute the values of
$v \left( j_{E} \right)=v\left( c_{4}^{3}(E)/\Delta_{E} \right)$, and find that if
the Kodaira type of $E$ is $III$, $IV$, $III^{*}$, $IV^{*}$ or $I_{0}^{*}$,
then $v \left( j_{E} \right) \geq 0$. For Kodaira type $I_{m}^{*}$ with $m \geq 1$,
if $p \geq 3$, then $v\left( \Delta_{E} \right)=6+m$ and $v \left( c_{4}(E) \right)=2$,
so $v(j_{E})=-m<0$. Therefore such $E$ have potential good reduction if and only
if $p=2$ and $v \left( j_{E} \right) \geq 0$.

We write $x(P)=x$ and $x(P')=x'$ for convenience.

\begin{lemma}
\label{lem:multiplicative red}
Suppose that $E/K$ is an elliptic curve having additive reduction and $P \in E(K)$
having singular reduction.

\noindent
{\rm (i)} There exists a finite extension, $L$, of $K$ such that $E$ has either
good or multiplicative reduction over $K$. So there exists a change of variables,
$x=u^{2}x'+r$ and $y=u^{3}y'+u^{2}sx'+t$, from $E$ to an elliptic curve $E'$ in
minimal Weierstrass form with $r,s,t,u \in R_{L}$. Let $P'$ be the image of $P$
under this change of variables.

\noindent
{\rm (ii)} Let $P'$ be as in part~{\rm (i)}. If $P'$ has singular reduction,
then $\ell_{P'}=2a_{P'}=v_{L}(\Delta_{E'})$, where $\ell_{P'}=-v(j(E'))$ and
$a_{P'}$ is the component of the N\'{e}ron model special fibre $(\cong \bbZ/\ell_{P'}\bbZ)$
containing $P'$.
\end{lemma}

\begin{proof}
(i) This is the semi-stable reduction theorem. See Proposition~VII.5.4
of \cite{Silv2}.

(ii) Since $P'$ has singular reduction, $E'$ has multiplicative reduction by
part~(i) and $E$ has Kodaira type $I_{m}^{*}$ with $m \geq 1$ (see, for example,
the passage below Th\'{e}or\`{e}me~2 on page~121 of \cite{Papa}) and
$v \left( j_{E} \right)=v\left(c_{4}^{3}(E)/\Delta_{E}\right)<0$. Here $c_{v}=2$
or $4$.

From Tableaux~I--III and
Th\'{e}or\`{e}me~3 on page~121 in \cite{Papa}, we find that
\begin{equation}
\label{eq:Delta1}
v \left( \Delta_{E} \right)=m+4+v\left(c_{4}(E)\right).
\end{equation}

From Lemma~\ref{lem:val-extensions}, we know that $v_{L}(u)=T(E)/2
=v_{L} \left( c_{4}(E) \right)/4$. Combining this with the expression for
$v_{L} \left( \Delta_{E'} \right)$ in Table~3.1 on page~45 of
\cite{Silv2}, it follows that
\begin{equation}
\label{eq:Delta2}
v_{L} \left( \Delta_{E'} \right)
= v_{L} \left( \Delta_{E} \right) - 3v_{L} \left(c_{4}(E) \right)
= (m+4)e_{L}-2v_{L} \left( c_{4}(E) \right).
\end{equation}

By Proposition~2.2(iv) of \cite{Stange} and $v_{L}(u)=v_{L} \left( c_{4}(E) \right)/4$,
\[
v_{L} \left( \psi_{2}\left( P' \right) \right)
=v_{L} \left( \psi_{2}\left( P \right) \right)-3v_{L}(u)
=v_{L} \left( \psi_{2}\left( P \right) \right)-\frac{3}{4}v_{L}\left(c_{4}(E)\right).
\]

By Lemma~5.1 of \cite{Silv4},
\[
a_{P'} = \min \left\{ v_{L} \left( \psi_{2}\left( P' \right) \right), v_{L}\left(\Delta_{E'}\right)/2 \right\}.
\]
Applying \eqref{eq:Delta1} and then \eqref{eq:Delta2}, we have
\[
\ell_{P'} = -v_{L}\left( j_{E} \right)=-v_{L}\left(c_{4}^{3}(E)/\Delta_{E}\right)
=(m+4)e_{L}-2v_{L}\left(c_{4}(E)\right)=v_{L} \left( \Delta_{E'} \right).
\]

We consider two cases in order to determine $a_{P'}$.

\vspace*{2.0mm}

(a) First, assume that $[2]P$ has non-singular reduction. We establish some
information about $n_{P'}$.

If $v_{L}\left(x\left(P'\right)\right) < 0$, then $n_{P'}=1$.
 
Next we consider $v_{L}\left(x\left(P'\right)\right) \geq 0$. We saw above that
$v_{L} \left( u^{2} \right)=v_{L} \left( c_{4}(E) \right)/2$. So $v_{L} \left( u^{2} \right)>0$,
by Proposition~VII.5.1(c) of \cite{Silv2}.
Since $v_{L}\left(x\left(P'\right)\right) \geq 0$ and $v_{L} \left( u^{2} \right) > 0$,
it follows that $v_{L}(x(P)-r) \geq v_{L} \left( u^{2} \right) > 0$. Thus, since
$v_{L}(x(P))>0$ by Lemma~\ref{lem:sing-pts}(i), we must also have $v_{L}(r)>0$.
Since $[2]P$ is non-singular, by Lemma~\ref{lem:sing-pts}(i) we have $v_{L} \left( x([2]P) \right) \leq 0$. As a result,
\[
v_{L} \left( x \left( [2]P' \right) \right)
= v_{L} \left( u^{-2} \left( x([2]P) - r \right) \right)
\leq -2v_{L}(u)<0,
\]
and hence $n_{P'} | 2$. 

With this information about $n_{P'}$, we proceed to determine $a_{P'}$.

Suppose that $v_{L} \left( \psi_{2}\left( P' \right) \right) < v_{L}\left(\Delta_{E'}\right)/2$.
Then $a_{P'}=v_{L} \left( \psi_{2}\left( P' \right) \right)$.
From $\ell_{P'}=v_{L} \left( \Delta_{E'} \right)$, we have $\ell_{P'}>2a_{P'}$.
Then $\widehat{a_{P'}}=a_{P'}$ and $\widehat{2a_{P'}}=2a_{P'}$, where $\widehat{x}$
denotes the least non-negative residue of $x$ modulo $\ell_{P'}$. Hence
\[
R_{2}\left(a_{P'}, \ell_{P'}\right)
=\frac{4a_{P'}\left(\ell_{P'}-a_{P'}\right)}{2\ell_{P'}}
 -\frac{2a_{P'} \left( \ell_{P'}-2a_{P'} \right)}{2\ell_{P'}}=a_{P'},
\]
where $R_{n}$ is as defined in Definition~\ref{def:Rn}.

From Lemma~\ref{lem:phi}(ii), we have $v_{L}\left(\phi_{2}\left(P'\right)\right)
=2R_{2}\left( a_{P'}, \ell_{P'} \right)$. So $v_{L}\left(\phi_{2}\left(P'\right)\right)
=2a_{P'}=v_{L} \left( \psi_{2}^{2}\left( P' \right) \right)$. Therefore
$v_{L} \left( x \left( [2]P' \right) \right)=0$, which contradicts $n_{P'}=1$ or $2$.
So $v_{L} \left( \psi_{2}\left( P' \right) \right) \geq v_{L}\left(\Delta_{E'}\right)/2$.
Hence $a_{P'}=v_{L} \left( \Delta_{E'} \right)/2$ and it follows that $\ell_{P'}=2a_{P'}$. 

\vspace*{2.0mm}

(b) Next, assume that $[2]P$ has singular reduction. From Lemma~\ref{lem:ImStar-vals},
this only happens when $m$ is odd and $v\left( \psi_{2}^{2}(P) \right)=m+3$.
Using this, Proposition~2.2(iv) of \cite{Stange} (again noting that our change
of variables is defined using Silverman's convention, which differs from Stange's
in \cite{Stange}) and $v_{L}(u)=v_{L} \left( c_{4}(E) \right)/4$,
we find that
\[
v_{L} \left( \psi_{2}\left( P' \right) \right) = v_{L} \left( \psi_{2}(P) \right)-3v_{L}(u)
=\frac{(m+3)e_{L}}{2}-\frac{3v_{L}\left(c_{4}(E)\right)}{4}.
\]

We find that $v \left( c_{4}(E) \right) \geq 2$ from Tableaux~I--V of \cite{Papa}.
Applying this and \eqref{eq:Delta2} to the preceding expression for
$v_{L} \left( \psi_{2}\left( P' \right) \right)$, we obtain
\begin{align*}
v_{L} \left( \psi_{2}\left( P' \right) \right)-\frac{v_{L} \left( \Delta_{E'} \right)}{2}
=\frac{v_{L} \left( c_{4}(E) \right)-2e_{L}}{4} \geq 0.
\end{align*}
So $v_{L} \left( \psi_{2}\left( P' \right) \right) \geq v_{L} \left( \Delta_{E'} \right)/2$.
Hence $a_{P'}=v_{L} \left( \Delta_{E'} \right)/2$ and it follows that $\ell_{P'}=2a_{P'}$.
\end{proof} 

\begin{lemma}
\label{lem:add-relations}
Suppose that $E/K$ is an elliptic curve having additive reduction and $P \in E(K)$
having singular reduction. Let $E'$, $L$, $P'$ and $r$ be as in
Lemma~{\rm \ref{lem:multiplicative red}(i)}.

\noindent
{\rm (i)} If $n \equiv 0 \bmod{m_{P}}$ and $P'$ has non-singular reduction, then 
\begin{equation}
\label{eq:kpn-1}
k_{v,n}(P)=v\left( \phi_{n}(P) \right) =\left\{
\begin{array}{ll}
v(x(P)- r) n^{2}             & \text{if $v_{L}\left(x\left(P'\right)\right) \leq 0$,} \vspace{1mm}\\ 
v \left( u^{2} \right) n^{2} & \text{if $v_{L}\left(x\left(P'\right)\right)>0$.}
\end{array}
\right. 
\end{equation}

\noindent
{\rm (ii)} If $n \equiv 0 \bmod{m_{P}}$ and $P'$ has singular reduction, then 
\begin{equation}
\label{eq:kpn-2}
k_{v,n}(P)=v\left( \phi_{n}(P) \right) =\frac{1}{4}(m+4) n^{2}.
\end{equation}

\noindent
{\rm (iii)} If $n \equiv \pm 1 \bmod{m_{P}}$ and $P'$ has non-singular reduction, then 
\begin{equation}
\label{eq:kpn-3}
k_{v,n}(P)=v\left( \psi_{n}^{2}(P) \right) =v(x(P)-r) (n^{2}-1).
\end{equation}

\noindent 
{\rm (iv)} If $n \equiv \pm 1 \bmod{m_{P}}$ and $P'$ has singular reduction, then
\begin{equation}
\label{eq:kpn-4}
k_{v,n}(P)=v\left( \psi_{n}^{2}(P) \right) =\frac{1}{4}(m+4) (n^{2}-1).
\end{equation}
\end{lemma}

\begin{proof}
We will first show that when $n \equiv 0 \bmod{m_{P}}$, we have $n_{P'}|n$.

If $v_{L}\left(x\left(P'\right)\right) < 0$, then $n_{P'}=1$, and so $n_{P'} \mid n$.
If $v_{L}\left(x\left(P'\right)\right) \geq 0$, then using the same argument as
in case~(a) of the proof of Lemma~\ref{lem:multiplicative red}(ii), we find that
$n_{P'}|n$.

\vspace*{2.0mm}
   
(i) Since $v \left( x\left( [n]P \right) \right) \leq 0$, we have
$v\left( \phi_{n}(P) \right) \leq v\left( \psi_{n}^{2}(P) \right)$ by
Lemma~\ref{lem:sing-pts}(i). Therefore $k_{v,n}(P)=v \left( \phi_{n}(P) \right)$.

We first assume that $v_{L}\left(x\left(P'\right)\right) \leq 0$. We can apply
Lemma~\ref{lem:phi}(i), and obtain $v_{L} \left( \phi_{n}\left(P'\right) \right)
= v_{L} \left( x\left(P'\right) \right) n^{2}$. So by \eqref{eq:lemma 13} in
Lemma~\ref{lem:val-extensions}(ii) with $T(E)=v_{L} \left(u^{2}\right)$,
\begin{align*}
v_{L}\left( \phi_{n}(P) \right)
&= \left( n^{2}-1 \right)v_{L}\left(u^{2}\right) + v_{L} \left( u^{2}\phi_{n}\left( P' \right) \right)
= \left(v_{L}\left(u^{2}\right)+v_{L}\left(x\left(P'\right)\right)\right)n^{2} \\
&= v_{L}\left( x(P)- r \right) n^{2}.
\end{align*}

Now assume that $v_{L}\left(x\left(P'\right)\right)> 0$. Applying Lemma~\ref{lem:phi}(i),
we obtain $v_{L} \left( \phi_{n}\left(P'\right) \right)=0$. As above, we have
\[
v_{L}\left( \phi_{n}(P) \right)
= \left( n^{2}-1 \right)v_{L}\left(u^{2}\right) + v_{L} \left( u^{2}\phi_{n}\left( P' \right) \right)
=v_{L}\left( u^{2} \right)n^{2}.
\]

\vspace*{2.0mm}

(ii) Since $P'$ has singular reduction, it follows that $E'$ has multiplicative reduction.
From Lemma~\ref{lem:multiplicative red}(ii), we know that $\ell_{P'}=2a_{P'}
=v_{L} \left( \Delta_{E'} \right)$.

Since $m_{P}$ is even for $I_{m}^{*}$, from our assumption that $n \equiv 0 \bmod{m_{P}}$,
we see that $n$ is even. So $na_{P'} \equiv 0 \bmod{\ell_{P'}}$. Therefore
\begin{equation}
\label{eq:rn-even}
R_{n} \left( a_{P'},\ell_{P'} \right) = \frac{n^{2}a_{P'} \left( 2a_{P'}-a_{P'} \right)}{4a_{P'}}
=\frac{a_{P'}n^{2}}{4}.
\end{equation}

Since $n \equiv 0 \bmod{m_{P}}$, we have $k_{v,n}(P)=v\left(\phi_{n}(P)\right)$.
Since $E'$ has multiplicative reduction and $n_{P'} | n $, from Lemma~\ref{lem:phi}(ii)
we have $v_{L} \left( \phi_{n}\left(P'\right) \right)=2R_{n}\left( a_{P'}, \ell_{P'} \right)$.
From \eqref{eq:lemma 13} in Lemma~\ref{lem:val-extensions}, followed by \eqref{eq:Delta2} and 
\eqref{eq:rn-even}, and noting that
$v_{L}(u)=v_{L}\left(c_{4}(E)\right)/4$, we obtain
\begin{align*}
v_{L} \left( \phi_{n}(P) \right)
&= v_{L}\left(u^{2}\right) \left( n^{2}-1 \right) + v_{L} \left( u^{2}\phi_{n}\left(P'\right) \right) \\
&= v_{L}\left(u^{2}\right)n^{2}+2R_{n}\left(a_{P'}, \ell_{P'}\right) \\
&= v_{L}\left(u^{2}\right)n^{2}+\frac{a_{P'}n^{2}}{2} \\
&= v_{L}\left(u^{2}\right)n^{2}+\frac{v_{L} \left( \Delta_{E'} \right)}{4} \\
&= \frac{1}{2} v_{L}\left(c_{4}(E)\right) n^{2} +\frac{n^{2}}{2} \left(\frac{(m+4)e_{L}}{2}-v_{L}\left(c_{4}(E)\right)\right) \\
&= \frac{1}{4}(m+4)n^{2}e_{L}.
\end{align*}

Hence
\[
k_{v,P}(P)=v \left( \phi_{n}(P) \right)=\frac{1}{4}(m+4)n^{2}.
\]

\vspace*{2.0mm}

If $n \equiv \pm 1 \bmod{m_{P}}$, then $[n]P$ has singular
reduction and $v \left( x\left( [n]P \right) \right)>0$,
by Lemma~\ref{lem:sing-pts}(i). Hence $v\left( \phi_{n}(P) \right)
> v\left( \psi_{n}^{2}(P) \right)$ and so
$k_{v,n}(P)=v \left( \psi_{n}^{2}(P) \right)$.

We can write $n=k \pm 1$ for some positive integer $k$ with $k \equiv 0 \bmod{m_{P}}$.

(iii) We break the proof into two parts, depending on
$v_{L} \left( x\left(P'\right) \right)$.

Assume first that $v_{L}\left(x\left(P'\right)\right) \geq 0$. By the same
argument as at the start of the proof, we have $n_{P'} \mid k$ and
since $n_{P'}>1$ (by our assumption that $v_{L}\left(x\left(P'\right)\right) \geq 0$),
we have $n_{P'} \nmid n$. So by Lemma~\ref{lem:non-sing} we have $v_{L} \left( \psi_{n} \left( P' \right) \right)=0$.
By using \eqref{eq:psi-extension} in Lemma~\ref{lem:val-extensions}
with $T(E)=v_{L} \left( u^{2} \right)$, we obtain
\[
v_{L} \left( \psi_{n}^{2}(P) \right)
= v_{L} \left( u^{2} \right) \left( n^{2}-1 \right)
= v_{L}(x(P)-r) \left( n^{2}-1 \right).
\]

Now assume that $v_{L}\left(x\left(P'\right)\right)<0$. Then $n_{P'}=1$ and $n_{P'} | n$.
From Lemma~\ref{lem:non-sing},
\[
v_{L}\left( \psi_{n}^{2}\left( P' \right) \right) = v_{L} \left(x\left( P' \right)\right) n^{2}+2S_{n} \left(P'\right).
\]

From the definition of $S_{n}\left(P'\right)$ in \eqref{eq:s-defn} along with
the conditions in \eqref{eq:param-conditions}, we have
$S_{n}\left(P'\right) \geq s_{P'}
=v_{L} \left(x\left( P' \right)/y\left( P' \right) \right)$. Noting that
$3v_{L} \left(x\left( P' \right)\right)=2v_{L} \left(y\left( P' \right)\right)$,
we obtain $s_{P'}= -v_{L}\left(x \left( P' \right)\right)/2$. Hence
\[
v_{L}\left( \psi_{n}^{2}\left( P' \right) \right) \geq v_{L} \left( x\left( P' \right) \right) \left( n^{2}-1 \right)
\]
and so, from \eqref{eq:psi-extension} in Lemma~\ref{lem:val-extensions}(i),
\[
v_{L} \left( \psi_{n}^{2}(P) \right)
\geq v_{L} \left( u^{2}\right) \left( n^{2}-1 \right)
     + v_{L} \left(x\left( P' \right)\right) \left( n^{2}-1 \right)
=    v_{L}(x(P)-r) \left( n^{2}-1 \right).
\]

We can write
\begin{equation}
\label{eq:psi-rel}
v_{L} \left( \psi_{n}^{2}(P) \right)=v_{L} \left( x(P)-r \right) \left( n^{2}-1 \right)+\alpha_{n},
\end{equation}
where $\alpha_{n} \geq 0$. Since $n_{P'}=1$, this argument works for any value
of $n$, not just $n \equiv \pm 1 \bmod{m_{p}}$. Hence \eqref{eq:psi-rel} holds
for all $n$ too.

We will now prove that, in fact, $v_{L} \left( \psi_{n}^{2}(P) \right)=v_{L}(x(P)-r) \left( n^{2}-1 \right)$.
Since $v_{L}(x([k]P)) \leq 0$ and $v_{L}(x(P))>0$,
we have $v_{L} \left( \phi_{k}(P) \right) \leq v_{L}\left(\psi_{k}^{2}(P)\right)
<v_{L}\left( x(P)\psi_{k}^{2}(P) \right)$, so
from the formula $\phi_{k}=x\psi_{k}^{2}-\psi_{k-1}\psi_{k+1}$ and applying
\eqref{eq:psi-rel} with $k-1$ and $k+1$ instead of $n$, we obtain
\begin{equation}
\label{eq:phi_k}
v_{L}\left(\phi_{k}(P)\right)=v_{L}\left(\psi_{k-1}(P)\psi_{k+1}(P)\right)
=v_{L} \left( x(P)-r \right) k^{2}+\alpha_{k-1}/2+\alpha_{k+1}/2.
\end{equation}

On the other hand, since $k \equiv 0 \bmod{m_{P}}$, from part~(i) of this
lemma, we have $v_{L}\left(\phi_{k}(P)\right)
=v_{L} \left( x(P)-r \right)k^{2}$. Therefore
$\alpha_{k-1}=\alpha_{k+1}=0$. 

\vspace*{2.0mm}

(iv) Here $E'$ has multiplicative reduction.
So by Lemma~\ref{lem:multiplicative red}(ii), we know that $a_{P'}=v_{L}(\Delta_{E'})/2$
and $\ell_{P'}=2a_{P'}$.

Since $m_{P}=2$ or $4$, we see that $n$ is odd. 
So $na_{P'} \equiv a_{P'} \bmod{2a_{P'}}$. From $\ell_{P'}=2a_{P'}$, we have
$\widehat{a_{P'}}=a_{P'}$ and $\widehat{na_{P'}}=a_{P'}$, where $\widehat{x}$
denotes the least non-negative residue of $x$ modulo $\ell_{P'}$. Hence
\begin{equation}
\label{eq:n-odd}
v_{L}\left(\psi_{n}\left(P'\right)\right)=R_{n}\left(a_{P'}, \ell_{P'}\right)
=\frac{a_{P'}(n^{2}-1)}{4}.
\end{equation}

By \eqref{eq:psi-extension} in Lemma~\ref{lem:val-extensions}, followed by \eqref{eq:Delta2} and 
\eqref{eq:n-odd}, and noting that $v_{L}(u)=v_{L}\left(c_{4}(E)\right)/4$,
we obtain
\begin{align*}
v_{L} \left( \psi_{n}^{2}(P) \right)
&= \frac{1}{2} v_{L}\left(c_{4}(E)\right) \left( n^{2}-1 \right)
   +\frac{1}{2} \left( n^{2}-1 \right) \left(\frac{(m+4)e_{L}}{2}-v_{L}\left(c_{4}(E)\right)\right) \\
&= \frac{1}{4}(m+4)e_{L} \left( n^{2}-1 \right).
\end{align*}

Hence 
\[
k_{v,n}(P)=v \left( \psi_{n}^{2}(P) \right)=\frac{1}{4}(m+4) \left( n^{2}-1 \right).
\]
\end{proof}


\begin{lemma}
\label{lem:k-p,n}
Suppose that $E/K$ is an elliptic curve having additive reduction and $P \in E(K)$
having singular reduction.

{\rm (i)} If $c_{v}=2$, then
\[
k_{v,n}(P)=
\left\{
\begin{array}{ll}
v \left( \psi_{3}(P) \right)n^{2}/4                  & \text{if $n \equiv 0 \bmod{2}$,} \\
v \left( \psi_{3}(P) \right)\left( n^{2}-1 \right)/4 & \text{if $n \equiv 1 \bmod{2}$.}
\end{array}
\right.
\]

{\rm (ii)} If $c_{v}=3$, then
\[
k_{v,n}(P)=
\left\{
\begin{array}{ll}
v \left(\psi_{2}^{2}(P)\right)n^{2}/3                  & \text{if $n \equiv 0 \bmod{3}$,}\vspace*{1.0mm}\\
v \left(\psi_{2}^{2}(P)\right)\left( n^{2}-1 \right)/3 & \text{if $n \equiv 1,2 \bmod{3}$.}
\end{array}
\right.
\]

{\rm (iii)} Assume that $c_{v}=4$.

If $[2]P$ has non-singular reduction, then
\[
k_{v,n}(P)=
\left\{
\begin{array}{ll}
	v\left( \psi_{3} (P)\right)n^{2}/4                  & \text{if $n \equiv 0 \bmod{2}$,} \\
	v\left( \psi_{3} (P)\right)\left( n^{2}-1 \right)/4 & \text{if $n \equiv 1 \bmod{2}$.}
\end{array}
\right.
\]

If $[2]P$ has singular reduction, then
\[
k_{v,n}(P)
=\left\{
\begin{array}{ll}
	v\left( \psi_{3}(P) \right)n^{2}/4                  & \text{if $n \equiv 0 \bmod{4}$,} \\
	v\left( \psi_{3}(P) \right)\left( n^{2}-1 \right)/4 & \text{if $n \equiv 1,3 \bmod{4}$,} \\
	v\left( \psi_{3}(P) \right) n^{2} /4-1              & \text{if $n \equiv 2 \bmod{4}$.}	
\end{array} \right.
\]
\end{lemma}

\begin{proof}
Let $E'$, $L$, $P'$ and the change of variables be as in Lemma~\ref{lem:multiplicative red}(i).

\vspace*{1.0mm}

(i) Since $c_{v}=2$, it follows that $E$ has Kodaira type $III$, $III^{*}$ or
$I_m^{*}$ and that $m_{P}=2$, since $P$ has singular reduction.

(i-a) We first assume that $P'$ has non-singular reduction and $v_{L}\left(x\left(P'\right)\right)>0$.
From \eqref{eq:kpn-1}
with $n=2$, we have $v \left(\phi_{2}(P) \right)=4 v \left( u^{2} \right)$.
From \eqref{eq:kpn-3} with $n=3$, we have $v \left(\psi_{3}(P) \right)=4 v \left(x(P) -r\right)$.
From Lemma~\ref{lem:pot-good-red-vals}(i) along with Lemma~\ref{lem:ImStar-vals},
we have $v \left(\phi_{2}(P) \right)=v \left(\psi_{3}(P) \right)$, so
$v \left(\psi_{3}(P) \right)=v \left( u^{2} \right)$ if $n \equiv 0 \bmod 2$ and
$v \left(\psi_{3}(P) \right)=v \left(x(P) -r\right)$ if $n \equiv 1 \bmod 2$.

Therefore, by parts~(i) and (iii) of Lemma~\ref{lem:add-relations},
$k_{v,n}(P)=v \left( u^{2} \right)n^{2}=v \left(\psi_{3}(P) \right)n^{2}/4$ if
$n \equiv 0 \bmod 2$ and $k_{v,n}(P)=v \left(x(P)-r\right) \left( n^{2}-1 \right)
=v \left(\psi_{3}(P) \right)\left( n^{2}-1 \right)/4$ if $n \equiv 1 \bmod 2$.

\vspace*{1.0mm}

(i-b) Next assume that $P'$ has non-singular reduction and $v_{L}\left(x\left(P'\right)\right) \leq 0$.
From \eqref{eq:kpn-3} with $n=3$, we have $v \left(\psi_{3}(P) \right)=4 v \left(x(P) -r\right)$.

Therefore, by parts~(i) and (iii) of Lemma~\ref{lem:add-relations},
$k_{v,n}(P)=v \left( x(P) -r \right)n^{2}=v \left(\psi_{3}(P) \right)n^{2}/4$ if
$n \equiv 0 \bmod 2$ and $k_{v,n}(P)=v \left(x(P)-r\right) \left( n^{2}-1 \right)
=v \left(\psi_{3}(P) \right)\left( n^{2}-1 \right)/4$ if $n \equiv 1 \bmod 2$.

\vspace*{1.0mm}

(i-c) Next assume that $P'$ has singular reduction.
This can only happen if $E$ has Kodaira type $I_{m}^{*}$, so by
Lemma~\ref{lem:ImStar-vals}, we have again $v \left(\psi_{3}(P) \right)=m+4$.
Therefore, by parts~(ii) and (iv) of Lemma~\ref{lem:add-relations}, we obtain
the desired result here too.

\vspace*{2.0mm}

(ii) Since $c_{v}=3$, $E$ has Kodaira type $IV$ or $IV^{*}$ and $m_{P}=3$, since
$P$ has singular reduction. This case happens when $E$ has potential good
reduction. So $P'$ has non-singular reduction.

From \eqref{eq:kpn-1} with $n=3$, we have either
$v \left(\phi_{3}(P) \right)=9 v\left( u^{2} \right)$ or
$v \left(\phi_{3}(P) \right)=9 v\left( x(P)-r \right)$. In either case, we
find by \eqref{eq:kpn-1} that
$v \left( \phi_{n}(P) \right)=v \left(\phi_{3}(P) \right) n^{2}/9$ holds for
$n \equiv 0 \bmod 3$.
By part~(ii) of Lemma~\ref{lem:pot-good-red-vals},
we have $v \left( \phi_{3}(P) \right)=3v \left( \psi_{2}^{2}(P) \right)$,
so $v \left( \phi_{n}(P) \right)=v \left(\psi_{2}^{2}(P) \right) n^{2}/3$.

From
\eqref{eq:kpn-3} with $n=2$, we have $v \left(\psi_{2}^{2}(P) \right)=3 v \left(x(P) -r\right)$.
So if $n \equiv \pm 1 \bmod 3$, then
$v(\psi_{n}(P))=v \left(\psi_{2}^{2}(P) \right) (n^{2}-1)/3$, using
\eqref{eq:kpn-3} again.

\vspace*{1.0mm}

(iii) Lastly, we consider $c_{v}=4$. This case happens only when the Kodaira
type is $I_{m}^{*}$.

(iii-a) Assume that $[2]P$ has non-singular reduction. Thus, $m_{P}=2$
and $v\left( \phi_{2}(P) \right) \leq v\left( \psi_{2}^{2}(P) \right)$ by
Lemma~\ref{lem:sing-pts}(i). Since $P$ has singular reduction, we have
$v(x(P)) \geq 1$ (again, by Lemma~\ref{lem:sing-pts}(i)).

If $n \equiv 0 \bmod{2}$, then \eqref{eq:kpn-1} or \eqref{eq:kpn-2} holds.
As shown in the proof of part~(i), $k_{v,n}=v\left( \psi_{3}(P) \right)n^{2}/4$.

Similarly, if $n \equiv 1 \bmod{2}$, then we can apply \eqref{eq:kpn-3} and \eqref{eq:kpn-4}
as in the proof of part~(i) to show that
$k_{v,n}=v\left( \psi_{3}(P) \right) \left( n^{2}-1 \right)/4$.

(iii-b) Assume that $[2]P$ has singular reduction. From Lemma~\ref{lem:ImStar-vals}(ii),
this case happens only when $m$ is odd. In this case, $m_{P}=4$.


If $n \equiv \pm 1 \bmod{4}$, then $[n]P$ is singular, so $k_{v,n}=v \left( \psi_{n}^{2}(P) \right)$.
By Lemma~\ref{lem:ImStar-vals}(i), $v\left(\psi_{3}(P) \right)=m+4$.
If $P'$ has singular reduction, then by \eqref{eq:kpn-4}, we obtain
$k_{v,n}=v\left( \psi_{3}(P) \right)\left( n^{2} - 1 \right)/4$.

If $P'$ has non-singular reduction, then by \eqref{eq:kpn-3}, we obtain
$k_{v,n}=v\left( x(P)-r \right)\left( n^{2} - 1 \right)$.
Applying \eqref{eq:kpn-3}, but with $n=3$, we find that
$8v\left( x(P)-r \right)=2v\left( \psi_{3}(P) \right)$, so the desired result
follows.

Assume that $n \equiv 0 \bmod{4}$. If $P'$ has singular reduction, then
\eqref{eq:kpn-2} applies. So, from Lemma~\ref{lem:ImStar-vals}(i)
we have $v\left( \psi_{3}(P) \right)=m+4$. Therefore $k_{v,n}=v\left( \psi_{3}(P) \right)n^{2}/4$.

If $P'$ has non-singular reduction,
then by \eqref{eq:kpn-1}, we obtain
$k_{v,n}=v\left( x(P)-r \right)n^{2}$
if $v_{L}\left(x\left(P'\right)\right) \leq 0$
and
$k_{v,n}=v\left( u^{2} \right)n^{2}$,
if $v_{L}\left(x\left(P'\right)\right)>0$.

We saw above that
$8v\left( x(P)-r \right)=2v\left( \psi_{3}(P) \right)$, so
$k_{v,n}=v\left( \psi_{3}(P) \right)n^{2}/4$, if
$v_{L}\left(x\left(P'\right)\right) \leq 0$.

Applying \eqref{eq:kpn-1} with $n=4$, we find that
$16v\left( u^{2} \right)=v\left( \phi_{4}(P) \right)$ if
$v_{L}\left(x\left(P'\right)\right)>0$. So
$k_{v,n}=v\left( \phi_{4}(P) \right)n^{2}/16$, if
$v_{L}\left(x\left(P'\right)\right)>0$.

Using \eqref{eq:phi_k} with $k=4$ (note that it is applicable since
$[4]P$ is non-singular and $P$ is singular), we have
$v \left( \phi_{4}(P) \right) = v \left( \psi_{3}(P)\psi_{5}(P) \right)$.
By \eqref{eq:kpn-3} and part~(iii) of this lemma with $n=5$ (note that we have
already proven part~(iii) for such $n$), we have
$v \left( \psi_{5}^{2}(P) \right)= 6v \left( \psi_{3}(P) \right)$, so
$v \left( \psi_{5}(P) \right)= 3v \left( \psi_{3}(P) \right)$. Hence
$v \left( \phi_{4}(P) \right) = 4v \left( \psi_{3}(P) \right)$, and our desired
result follows.

Assume that $n \equiv 2 \bmod{4}$.
Since $[n]P$ has singular reduction and $[2n]P$ has non-singular reduction,
from Lemma~\ref{lem:ImStar-vals}(i), it must be that $v \left( x\left( [n]P \right) \right)=1$.
Therefore $v \left( \psi_{n}^{2}(P) \right) =v\left(\phi_{n}(P)\right)-1$.
Lemma~\ref{lem:ImStar-vals}(i) tells us that $v(x(P)) \geq 2$,
so $v\left( \phi_{n}(P) \right)<v\left( x(P)\psi_{n}^{2}(P) \right)$.
From the formula $\phi_{n}=x\psi_{n}^{2}-\psi_{n-1}\psi_{n+1}$ along with the
expressions just found when $n \equiv \pm 1 \bmod{4}$, we obtain
\begin{align*}
v \left( \phi_{n} \left( P \right) \right)
&= v\left(\psi_{n-1}(P)\psi_{n+1}(P)\right) \\
&= \frac{1}{4}v \left( \psi_{3}(P) \right)\left( \left( n-1 \right)^{2}-1 \right)
   +\frac{1}{4}v\left( \psi_{3}(P) \right)\left( \left( n+1 \right)^{2}-1\right) \\
&= \frac{1}{4}v\left( \psi_{3}(P) \right)n^{2}.
\end{align*}
\end{proof}

We now show how Theorem~\ref{thm:main} follows when $E$ has additive reduction.

If $E$ has Kodaira type $IV$ or $IV^{*}$, then 
\begin{equation}
\label{eq:val-psi 2}
\frac{1}{3}v\left(\psi_{2}^{2}(P) \right)=-2\frac{[K:\bbQ_{p}]}{\log |k|}\lambda_{v}(P),
\end{equation}
by Proposition~5(d) in \cite{CPS} (noting that our definition of $\lambda_{v}$,
following \cite{Silv4}, is one-half that in \cite{CPS}).

If $E$ has Kodaira type $III$, $III^{*}$, $I_{0}^{*}$ or $I_{m}^{*}$, then
\begin{equation}
\label{eq:val-psi 3}
\frac{1}{4}v\left(\psi_{3}(P) \right)=-2\frac{[K:\bbQ_{p}]}{\log |k|}\lambda_{v}(P),
\end{equation}
by Proposition~5(e) in \cite{CPS}.

Applying \eqref{eq:val-psi 2} and \eqref{eq:val-psi 3} with Table~2 of \cite{CPS}
to Lemma~\ref{lem:k-p,n}, we complete the proof of Theorem~\ref{thm:main}. For
$I_{2m}^{*}$ and $c_{v}=4$, we use $v\left(\phi_{2}(P) \right)=v\left(\psi_{3}(P) \right)$
from Lemma~\ref{lem:ImStar-vals} to get an expression in terms of $v\left(\phi_{2}(P) \right)$
rather than $v\left(\psi_{3}(P) \right)$.

\section*{Acknowledgements}

The authors greatly appreciate the suggestion from Federico Pellarin that we
investigate possible connections between the quantity, $k_{v,n}$, and local
heights. This led us to the current formulation of our Theorem~\ref{thm:main}.
We also thank the referee for their careful reading of our manuscript, helpful
suggestions for improving it and interesting remarks that broadened our understanding.

\bibliographystyle{amsplain}

\end{document}